# Real-Time Bi-directional Electric Vehicle Charging Control with Distribution Grid Implementation

Yingqi Xiong*, Behnam Khaki†, Chi-cheng Chu‡, Rajit Gadh§
Department of Mechanical and Aerospace Engineering
University of California, Los Angeles
Los Angele, USA
Email: *yxb936@g.ucla.edu, †behnamkhaki@g.ucla.edu, ‡peterchu@g.ucla.edu, §gadh@ucla.edu

*Abstract*—As electric vehicle (EV) adoption is growing year after year, there is no doubt that EVs will occupy a significant portion of transporting vehicle in the near future. Although EVs have benefits for environment, large amount of un-coordinated EV charging will affect the power grid and degrade power quality. To alleviate negative effects of EV charging load and turn them to opportunities, a decentralized real-time control algorithm is developed in this paper to provide optimal scheduling of EV bi-directional charging. To evaluate the performance of the proposed algorithm, numerical simulation is performed based on real-world EV user data, and power flow analysis is carried out to show how the proposed algorithm improve power grid steady state operation. . The results show that the implementation of proposed algorithm can effectively coordinate bi-directional charging by 30% peak load shaving, more than 2% of voltage drop reduction, and 40% transmission line current decrease.

## I. INTRODUCTION

As what we have seen in recent years, the fast development of Electric Vehicle (EV) technology and the stable and continuous growth of EV adoption make it obvious that EVs play an important role in the transportation system in the very near future. Economic studies forecast that EV sales will hit 41 million by the year of 2040, representing 35% of new light duty vehicle sales [1]. While charging station becomes more and more accessible with the amount of installation increases rapidly, bi-directional charging stations can discharge the electric energy in the EV battery back into the power grid, vehicle-to-grid (V2G), and provide ancillary service [2]. According to the SAE CCS and IEC15118 standards, the power limit of fast charging is up to 120kW, which can be equal to the total load of 30 to 40 households [3]. The gigantic energy consumption of EV chargers and V2G put challenges but also opportunities over the modern power grid. Although un-coordinated charging activities can cause voltage fluctuation, and power loss increase, especially when the number of charging sessions exceeds a certain level [4], coordinated smart charging can not only resolve the impact to the power grid, but also turn EV charging demand to an opportunity to flatten the profile of the system [5]. Moreover, charging stations with SAE CCS or CHAdeMO [6] protocol have V2G capability, which transforms EVs from simple dumb loads to distributed energy resources (DER), providing grid services such as voltage regulation and peak load shaving [7]. With the fast increasing EV possession rate, EVs will become a great asset of power grid as DERs in the near future, under the condition of being coordinated by well-designed control strategy. However, controlling thousands or millions of EV charging load is challenging, which comes into the following aspects: 1) EV charging scheduling is heavily constrained by user availability; 2) Huge computational resources are required for optimal charging scheduling of large number of EVs, especially in real-time; 3) The impact of smart bi-directional charging at transmission level is still a question.

Recent years we have seen a lot of research works proposing optimal controllers for coordinated charging of EVs and prediction methods for EV user behavior. Gan et al. [8] used a distributed protocol to optimally allocate EV charging time and energy. Wang et al. developed an EV charging scheduling algorithm with minimum cost as its objective [9]. Authors in [10] analyzed the influence of bi-directional EV charging to a building and provided related control strategy. Kernel density estimation method is used in [11] to deal with the uncertainty of EV user behavior based on real-world data. Wang et al. [12] utilized sample average approximation to predict EV schedule and performed a simulation with 10 EVs within a distribution grid. In [13] the author targeted the charging demand of EV in commercial parking lots and determine charging strategy using two-stage approximate dynamic programming. Charging coordination for Valley filling and cost reduction were studied in [14] and found out global optimum is not possible. However it is hard to find one among them which provides a comprehensive solution addressing all the aforementioned challenges based on real-world implementation requirements.

Accordingly, a decentralized control paradigm is introduced in this paper and applied to a power grid in order to investigate the influence of large scale EV charging load on the electrical grid. In our proposed decentralized algorithm, the computational burden of optimal charging scheduling is

This document was prepared as a result of work in part by grants from the California Energy Commission EPC-14-056 fund (Demonstration of PEV Smart Charging and Storage Supporting Grid Objectives Project). It does not necessarily represent the views of the Energy Commission, its employees, or the State of California. Neither the Commission, the State of California, nor the Commission's employees, contractors, nor subcontractors makes any warranty, express or implied, or assumes any legal liability for the information in this document; nor does any party represent that the use of this information will not infringe upon privately owned rights. This document has not been approved or disapproved by the Commission, nor has the Commission passed upon the accuracy of the information in this document.

distributed among the networked charging station equipped with embedded controller, which significantly reduces the computational resources required at control center compared to traditional centralized control. Western System Coordinating Council (WSCC) 9-bus model is used as the simulation test case, where the grid load profiles are generated by scaling up real world microgrid load profile collected by Cornell University [15], and they are assigned to the three PQ (load) buses of the grid model. EV user charging record data obtained from UCLA SMERC smart charging infrastructure [16] are used to model the EV charging load. The proposed decentralized charging control algorithm is performed over all EVs connected to the 9-bus system to flatten the total load profile in real time. Load flow analysis is carried out for worst cases before and after the implementation of proposed algorithm to evaluate the effectiveness of the proposed algorithm in improving voltage profile and reducing line current and generation cost. Specifically, the worst cases are defined as when large amount of EVs are plugged-in during peak hours while the grid baseload is also high. The contribution of this paper can be summarized as: 1) developing and implementing a real time decentralized optimal EV charging algorithm with real world data; 2) carrying out a large scale simulation to evaluate the impact of large number bi-directional charging activities to the power grid; 3) running load flow analysis to evaluate the performance of proposed algorithm in worst cases scenarios.

The rest of this paper is organized as follows: Section II discusses the smart charging infrastructure and configuration of the grid simulation. Section III explains the load flow analysis and real-time decentralized charging control algorithm. Section IV covers the numerical simulation results to evaluate the proposed algorithm performance, section V concludes the paper.

## II. SYSTEM OVERVIEW

### A. Smart Charging Infrastructure

SMERC EV charging system is used as test bed in this paper. This infrastructure has designed and implemented by UCLA Smart Grid Energy Research Center (SMERC) since 2013 and keeping evolving [4] [16]. It consists of 3 layers of components: control center, communication network and end devices such as charging stations and solar panels. There are more than 200 smart charging stations installed in UCLA campus, City of Santa Monica, and many other places around the Los Angeles area. All charging stations have a controller with embedded system and 3G network to send/receive data to/from control center. EV drivers use mobile app designed by SMERC to communicate with charging station and submit charging request. All charging activities are recorded and stored in the control center database. Prediction of EV user behavior and total EV load profile is made based on this data which has been collected for 4 years.

### B. Simulation Configuration

To investigate the impact of large scale EV charging activities on the distribution grid, scaled-up EV load profiles are generated based on the collected EV charging historical data. EV loads are then applied to WSCC 9-bus system model for load flow analysis and algorithm performance evaluation.

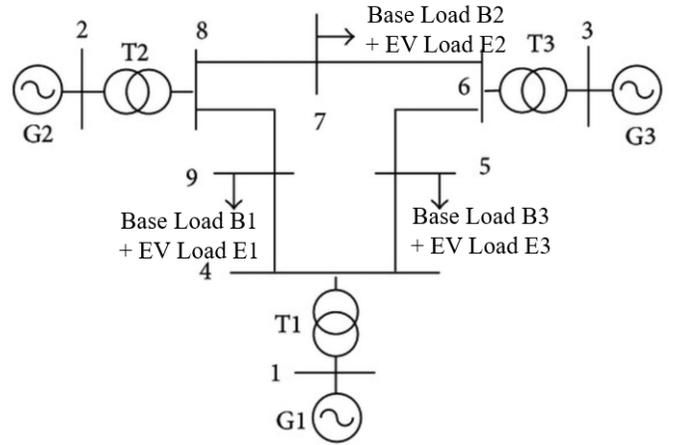

Figure 1. Single-line diagram of the distribution grid model

Three different base load and EV load profiles are applied to PQ buses 5, 7 and 9, as shown in Fig. 1. Voltage level at of these load buses is 230kV, and each of them represents a substation for a city region. Base load for each load bus is the total power consumption for the city without EVs, and EV load is the total power demand from all EV charging activities in that city. Three generators are connected to the grid through buses 1, 2, and 3 where bus#1 is the swing bus.

## III. LOAD FLOW AND CHARGING CONTROL

This section provides a brief description and review of power equation in load flow analysis. Formulation of the proposed decentralized EV charging algorithm is then explained in detail.

### A. Load Flow Analysis

Load flow analysis is conducted before and after the implementation of proposed algorithm to evaluate its contribution to the power quality. The relation between the current injected to each bus and the voltage of each bus of the grid model described in section II is shown by:

$$\begin{bmatrix} \bar{I}_1 \\ \vdots \\ \bar{I}_M \end{bmatrix} = \begin{bmatrix} \bar{Y}_{11} & \cdots & \bar{Y}_{1M} \\ \vdots & \ddots & \vdots \\ \bar{Y}_{M1} & \cdots & \bar{Y}_{MM} \end{bmatrix} \quad (1)$$

Or

$$\bar{I}_i = \sum_{j \in I} \bar{Y}_{ij} \bar{V}_j \quad (2)$$

Where $I_i$ is the injected current at bus $i$, and $V_j$ is the voltage at bus $j$. $Y_{ij}$ is the admittance of the line connecting bus $i$ to bus $j$, and $I = \{j | \bar{Y}_{ij} \neq 0\}$. The complex power injected to bus $i$ is obtained by:

$$\bar{S}_i = P_i + jQ_i = \bar{V}_i \bar{I}^*_i = \bar{V}_i \sum_{j \in I} \bar{Y}_{ij} \bar{V}_j \quad (3)$$

Where P and Q are active power and reactive power, respectively. The bus voltage and line admittance can be written as:

$$\bar{V}_i = V_i e^{j\theta_i} \quad (4)$$

$$\bar{Y}_{ij} = Y_{ij} e^{j\alpha_{ij}} \quad (5)$$

By replacing (4) and (5) in (3). Basic load flow equations are obtained.

$$\begin{cases} P_i = V_i \sum_{j \in I} Y_{ij} \cos(\theta_i - \theta_j - \alpha_{ij}) \\ Q_i = V_i \sum_{j \in I} V_j Y_{ij} \sin(\theta_i - \theta_j - \alpha_{ij}) \end{cases} \quad (6)$$

Newton-Raphson method is used to solve the above nodal power equations and determine the voltage at each load bus. By comparing the voltage drop during the time interval where large amount of EV charging activities taking places and the situation after charging control algorithm is applied, the impact of EV charging to the power grid and the effectiveness of scheduling algorithm can be evaluated.

### B. Decentralized EV Charging Control Algorithm

#### 1) Problem Overview

In this paper we assume that the utility or grid operator can always predict the day-ahead load profile in the region where their services are provided. Control center of the smart EV charging infrastructure retrieves the predicted base load as input data for the charging control algorithm. EVs are treated as distributed energy resources in the scope of this paper, which means the charging activities are bi-directional. EVs can either absorb energy from the power grid to feed their battery or discharge the energy in their battery to feed the power grid. The decision of charging or discharging is determined by proposed charging control algorithm based on the condition of current grid load and the EV driver demand. Control center learns EV user behavior from historical charging data to make prediction on parking schedule and possible state of charging (SoC) request, which are defined as the constraints in the optimal charging control algorithm. The decentralized control paradigm works as follows: 1) control center initializes and broadcasts control signal to all charging stations in its network; 2) each charging station updates its charging profile by performing an optimization only regarding to the EVs plugged in; 3) control center collects updated charging profiles from all charging stations, updates its control signal and broadcast it out again; 4) control algorithm is run in real-time with preset time interval and updated information. The goal of the proposed algorithm is to optimally schedule the bi-directional charging activities such that the load profile can be flattened, and load flow and voltage drop are kept in allowable range while all EV drivers' demands are fulfilled.

#### 2) Objective and Constraints

The prediction of baseload is denoted as $B^\tau(t)$ with $t \in T$ and $\tau \in H$. $T = [1, 2, \ldots T]$ is the time slot set where the algorithm will be performed, and $t$ is the specific time slot. $H$ is the time horizon where the charging control is moving forward. We assume that there are $N$ EVs supplied by the smart charging system. Each EV has a charging and discharging rate $p_n^\tau(t)$ with $n \in N = [1, 2, \ldots N]$ at time slot $t$ and time horizon $\tau$. The aim of proposed algorithm is to flatten the total load profile:

$$L^\tau(t) = \left( B^\tau(t) + \sum_{n=1}^N p_n^\tau(t) \right)^2 \quad (7)$$

We use $\overline{p_n}$ to denote maximum charging rate and $\overline{d_n}$ maximum V2G capacity of EV $n$. V2G discharging rate $\overline{d_n}$ is a negative value since the power flow is reversed. Thus we have:

$$\overline{d_n} \leq p_n^\tau(t) \leq \overline{p_n} \quad (8)$$

In our proposed control scheme, $\overline{p_n}$ and $\overline{d_n}$ are regarded as upper bound and lower bound on the controllable charging rate. During the time between an EV is plugged in the charging station, $t_{start}^{pred}$, and the time when it leaves, $t_{end}^{pred}$, charging rate is a continuous value between upper and lower bound, and it will be 0 otherwise, which means the specific EV is offline. EV charging rate and availability constraints are described in (9) and (10):

$$\overline{p_n} = \begin{cases} p_n^{\max}, t \in \left[ t_{start}^{pred}, t_{end}^{pred} \right] \\ 0, \quad otherwise \end{cases} \quad (9)$$

$$\overline{d_n} = \begin{cases} d_n^{\max}, t \in \left[ t_{start}^{pred}, t_{end}^{pred} \right] \\ 0, \quad otherwise \end{cases} \quad (10)$$

Energy consumptions of each EV are predicted from historical data. The predicted value $E_n^{pred}$ represents energy demand of EV user $n$ during his or her charging sessions. The summation of the products of EV charging rate and the corresponding time interval equals to the user energy demand:

$$\sum_T p_n(t) \Delta T = E_n^{pred} \quad (11)$$

Where $\Delta T$ is constant since each time slot is a fraction of the time horizon.

The control signal from control center is the first order derivative of total load profile which we intent to flatten, multiplied by a tuning parameter [7]:

$$c^i(t) = \frac{1}{\lambda N} \left( B^\tau(t) + \sum_{n=1}^N p_n^\tau(t) \right) \quad (12)$$

Where $i$ represents the iteration number in the algorithm, and $\lambda$ denotes the tuning parameter. The algorithm converges when the residual reaches convergence criteria, i.e. $\left\| c^{i+1} - c^i \right\| \leq \varepsilon$, where ε is a small positive real number.

Subsequently, by receiving the control signal, each charging station performs an optimization taking into consideration only the plugged-in EVs. The objective cost function is expressed as:

$$\sum_{t=1}^T c^i(t) p_n^{\tau^{i+1}}(t) + \frac{1}{2} \left\| p_n^{\tau^{i+1}}(t) - p_n^{\tau^i}(t) \right\|^2 \quad (13)$$

By minimizing the cost function, charging station will update EV charging profile $p_n^{\tau^{i+1}}(t)$ and report it to the control center.

#### 3) Real-time Implementation

The distributed algorithm to optimize EV charging energy allocation and scheduling is presented as following:

**Algorithm : Real Time Decentralized Bi-directional Charging**

Initialize all EV charging profiles $p_n^0(t) = 0$, $n = 1,2,\ldots N$

Pick control parameter $\lambda$

Initialize control signal $c^0(t) = \dfrac{1}{\lambda N}\left(B(t) + \sum_N p_n^0(t)\right)$

**While** $\tau < H$
    Updates prediction of baseload data
    Updates prediction of EV user behavior and energy demand
    **While** $\left\|c^{i+1} - c^i\right\| \geq \varepsilon$ :
        **For** each $n$ EVSE in the network, $n = 1,2,\ldots,N$
            Minimize (13)
            Subject to (8), (9), (10), (11)
        **End**
    $\tau = \tau + 1$
**End**

In the above algorithm, the control time horizon can be set to 24 hours, representing a day. The algorithm will update its control strategy every hour, i.e. $\tau \in [1,2,\ldots 24]$. When the algorithm converges, the total load profile defined in (5) will be effectively flattened. Fig. 2 shows the schematic of the proposed algorithm running in real time. Each charging station performs optimization by its own embedded controller, which greatly lowers the computational burden of the control center.

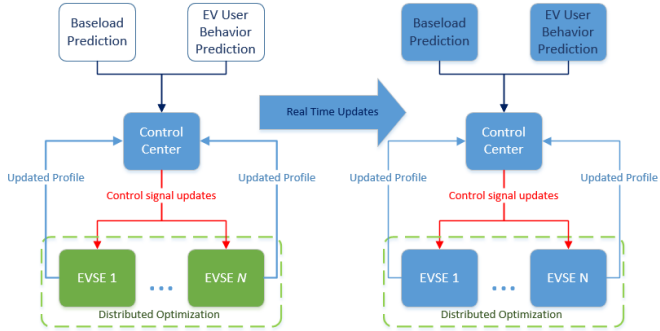

Figure 2. Schematic of proposed charging control algorithm

## IV. NUMERICAL SIMULATION

This section demonstrates and discusses the results obtained from large scale simulation of proposed decentralized charging control algorithm to show its capability in optimally scheduling EV charging activities and keeping power quality indices within the desired range.

### A. Test Case

WSCC 9-bus system is used as the test case for simulation in this paper. Hourly microgrid power consumption historical data is retrieved from Cornell University energy management portal [15] and scaled up to the level of a small city with around 50,000 population, which is comparable to the cities in Los Angeles area such as Culver City and Beverly Hills. The obtained load profiles are applied to the 3 load buses in the test case, namely bus#5, bus#7 and bus#9 as the base load. Prediction of EV user daily travelling schedule and energy demand is carried out based on the historical data in UCLA SMERC smart charging infrastructure collected over past 4 years, and it is utilized as constraints in the proposed algorithm. The number of EVs assigned to each load bus is 15,000, consisting of 35% light duty vehicles market shares as predicted for the year of 2040 in [1]. The simulation is carried out for different scenarios including 1) the worst case, where almost all EVs arrive and start charging at full capacity during peak demand hours; 2) the optimal case, where proposed algorithm shifts EV load and schedules V2G discharging sessions based on grid load condition and EV user behavior.

### B. Result and Analysis

In the condition that there is not any control or coordination over the EV charging stations, and all charging activities start immediately after EV is plugged in, total load can be too high for the power grid to stay in stable status.

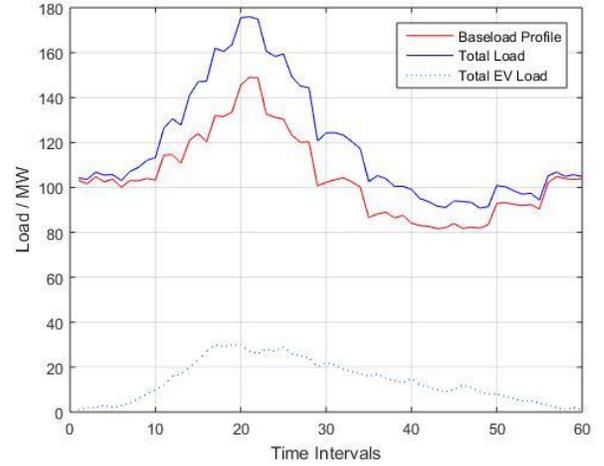

Figure 3. EV load peak overlap with grid peak demand

It can be seen in Fig. 3 that when no charging control is implemented, these is the potential that EV charging peak overlaps with gird peak demand (time interval~20), which is definitely not favorable to the power quality and power grid stability.

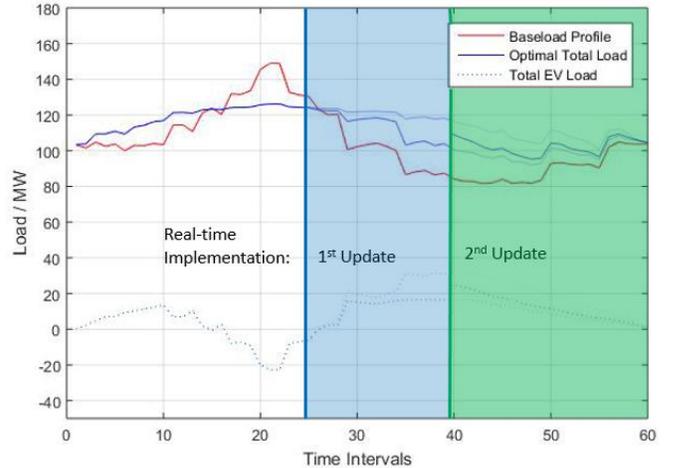

Figure 4. Optimal load profile with real time update

Fig. 4 shows the optimal total load profile, which is a combination of grid base load and total EV charging load, and it is much flatter and smoother than the base load is. By comparing Fig. 4 with un-coordinating charging scenario in Fig. 3, it can be seen that the proposed charging algorithm arranges as much V2G discharging sessions as during grid peak hours, and then coordinates full capacity charging activity over its control network during low load demand time intervals from 35 to 55. Another look into the total load profile reveals that the peak load after algorithm implementation is 125MW, while it is 180MW in the worst case without coordination. Real-time implementation is demonstrated in the shaded areas of Fig. 5, control center updates the EV user behavior prediction at time slot 25, according to current EV number and demand in the network, resulting in a new charging schedule and load profile. At time slot 40, second updates is made. New load profiles after updates are deeper blue color curves while previous profiles are shown in lighter color.

Load flow analysis results are shown in Table I, II and III. It can be seen that the total line current is reduced by 40% after charging control is implemented. The power drawn from swing bus is significantly decreased which results in generation cost reduction. Moreover, we observed 2% improvement of bus voltage profile over the grid while unacceptable voltage drops are seen at some of buses, such as bus#5 and bus#9, before EV charging coordination. Therefore, it can be said that EV charging load coordination for load flattening yields overall system operation improvement.

TABLE I. LINE CURRENT BEFORE / AFTER COORDINATION

| Line in the Grid | | Pre-control Line Current (A) | Post-control Line Current (A) |
|---|---|---|---|
| From | to | | |
| Bus 4 | Bus 9 | 420.8 | 130.7 |
| Bus 4 | Bus 5 | 420.4 | 153.7 |
| Bus 8 | Bus 7 | 358.9 | 232.8 |
| Bus 8 | Bus 9 | 58.7 | 167.7 |
| Bus 6 | Bus 7 | 185.0 | 75.8 |
| Bus 6 | Bus 5 | 41.4 | 132.5 |

TABLE II. POWER GENERATION BEFORE / AFTER COORDINATION

| Generator No. | Pre-control Generation Complex Power (MW) | Post-control Generation Complex Power (MW) |
|---|---|---|
| G1 (Swing) | 329.66 + j112.01 | 105 + j48.75 |
| G2 | 163 + j39.09 | 163 + j11.35 |
| G3 | 85 + j25.03 | 85 – j2.25 |

TABLE III. BUS VOLTAGE PROFILE BEFORE / AFTER COORDINATION

| Bus # | 4 | 5 | 6 | 7 | 8 | 9 |
|---|---|---|---|---|---|---|
| Before | 0.995 | 0.955 | 1.012 | 0.987 | 1.006 | 0.954 |
| After | 1.015 | 0.987 | 1.027 | 1.015 | 1.028 | 0.980 |

## V. CONCLUSION

In this paper, a real-time optimal decentralized EV charging algorithm is proposed and implemented. A large scale grid simulation is conducted with real-world load and EV user data to evaluate the algorithm performance. To evaluate the performance of the proposed algorithm from electrical grid point of view, power flow analysis is carried out for a test case with and without EV charging coordination. A 30% peak load shaving and more than 2% of voltage drop reduction are observed when EV charging load is controlled by the proposed algorithm. Also, Total line current of the power grid is reduced by 40% which significantly decreases the power loss.